\documentclass[11pt,a4paper]{amsart}

\usepackage{amssymb}
\usepackage{latexsym}
\usepackage{exscale}

\usepackage{mathrsfs}

\usepackage[colorlinks=true, pdfstartview=FitV, linkcolor=blue,
citecolor=blue, urlcolor=blue]{hyperref}

\newtheorem{theorem}{Theorem}[section]

\theoremstyle{definition}

\theoremstyle{remark}
\newtheorem{remark}[theorem]{Remark}

\numberwithin{equation}{section}

\renewcommand{\emptyset}{\mbox{\textup{\O}}}

\newcommand{\re}{\mathbb R}

\newcommand{\Z}{\mathbb Z}

\newcommand{\D}{\mathscr{D}}
\renewcommand{\S}{\mathcal{S}}

\DeclareMathOperator{\supp}{supp}

\newcommand{\Q}{\mathcal{Q}}

\def\Xint#1{\mathchoice
  {\XXint\displaystyle\textstyle{#1}}%
  {\XXint\textstyle\scriptstyle{#1}}%
  {\XXint\scriptstyle\scriptscriptstyle{#1}}%
  {\XXint\scriptscriptstyle\scriptscriptstyle{#1}}%
  \!\int}
\def\XXint#1#2#3{{\setbox0=\hbox{$#1{#2#3}{\int}$}
    \vcenter{\hbox{$#2#3$}}\kern-.5\wd0}}

\def\avgint{\Xint-}

\calclayout
\allowdisplaybreaks

\begin{document}

\title{A note on the off-diagonal Muckenhoupt-Wheeden conjecture}

\author{David Cruz-Uribe, SFO}
\address{David Cruz-Uribe, SFO\\
Dept. of Mathematics \\ Trinity College \\
Hartford, CT 06106-3100, USA} \email{david.cruzuribe@trincoll.edu}

\author{Jos\'e Mar{\'\i}a Martell}
\address{Jos\'e Mar{\'\i}a Martell
\\
Instituto de Ciencias Matem\'aticas CSIC-UAM-UC3M-UCM
\\
Consejo Superior de Investigaciones Cient{\'\i}ficas
\\
C/ Nicol\'as Cabrera, 13-15
\\
E-28049 Madrid, Spain} \email{chema.martell@icmat.es}

\author{Carlos P\'erez}
\address{Carlos P\'erez
\\
Departamento de An\'alisis Matem\'atico, Facultad de Mate\-m\'aticas\\
Universidad de Sevilla, 41012 Sevilla, Spain}
\email{carlosperez@us.es}

\subjclass{42B20, 42B25} \keywords{Haar shift operators, Calder\'on-Zygmund operators, two-weight inequalities, testing conditions.}

\thanks{The first author was supported by a grant from the Faculty
  Research Committee and the Stewart-Dorwart Faculty Development Fund
  at Trinity College; the first and third authors are supported by
  grant MTM2009-08934 from the Spanish Ministry of Science and
  Innovation; the second author is supported by grant Grant MTM2010-16518
  from the same institution}

\begin{abstract}
We obtain the off-diagonal Muckenhoupt-Wheeden conjecture for Calder\'on-Zygmund operators. Namely, given $1<p<q<\infty$ and a pair of weights $(u,v)$, if the Hardy-Littlewood maximal function satisfies the following two weight inequalities:
$$
M : L^p(v) \rightarrow L^q(u) \quad \text{ and } \quad
M: L^{q'}(u^{1-q'}) \rightarrow L^{p'}(v^{1-p'}),
$$
then any Calder\'on-Zygmund operator $T$ and its associated truncated maximal operator $T_\star$ are bounded from $L^p(v)$ to $L^q(u)$. Additionally, assuming only the second estimate for $M$ then $T$ and $T_\star$ map continuously $L^p(v)$ into $L^{q,\infty}(u)$. We also consider the case of generalized Haar shift operators and show that their off-diagonal two weight estimates are governed by the corresponding estimates for the dyadic Hardy-Littlewood maximal function.
\end{abstract}

\date{\today}

\maketitle

\section{Introduction and Main results}

In the 1970s, Muckenhoupt and Wheeden conjectured that given $p$, $1<p<\infty$, a sufficient
condition for the Hilbert transform to satisfy the two weight norm
inequality
\[ H : L^p(v) \rightarrow L^p(u) \]
is that the Hardy-Littlewood maximal operator satisfy the pair of norm
inequalities
\begin{gather*}
M : L^p(v) \rightarrow L^p(u), \\
M : L^{p'}(u^{1-p'}) \rightarrow L^{p'}(v^{1-p'}).
\end{gather*}
Moreover, they conjectured that the Hilbert transform satisfies the
weak-type inequality
\[ H : L^p(v) \rightarrow L^{p,\infty}(u) \]
provided that the maximal operator satisfies the second ``dual''
inequality.
Both of these conjectures  readily extend to all Calder\'on-Zygmund operators (see the definition below).  Very recently, both
conjectures  were disproved:  the strong-type inequality by Reguera
and Scurry~\cite{reguera-scurry} and the weak-type inequality by the
first author, Reznikov and Volberg~\cite{cruz-reznikov-volberg}.

\begin{remark}A special case of these conjectures, involving the  $A_p$ bump
conditions, has been considered by several authors:
see~\cite{MR2351373,MR2628851,cruz-martell-perezBook,DCU-JMM-CP2012,cruz-reznikov-volberg,lernerP2011-Dya}.
\end{remark}

In this note we prove the somewhat surprising fact that the
Muckenhoupt-Wheeden conjectures are true for off-diagonal inequalities.
Our main result is Theorem~\ref{thm:main-thm:CZO} below.  We also
prove an analogous result for the Haar shift operators (the so-called
dyadic Calder\'on-Zygmund operators) with the Hardy-Littlewood maximal operator
replaced by the dyadic maximal operator:  see
Theorem~\ref{thm:main-thm:dya} below.

To state our results we
first give some preliminary definitions. By weights we will always mean non-negative, measurable
functions.  Given a pair of weights $(u,v)$, hereafter we will assume
that $u>0$ on a set of positive measure and $u<\infty$ a.e., and $v>0$
a.e. and $v<\infty$ on a set of positive measure. We will also use the
standard notation $0\cdot\infty=0$.

\subsection*{Calder\'on-Zygmund operators}
A Calder\'on-Zygmund operator $T$ is a linear operator that is bounded
on $L^2(\re^n)$ and
$$
Tf(x)
=
\int_{\re^n} K(x,y)f(y)dy,
\qquad
f\in L^\infty_c(\re^n),
\quad
x\notin\supp f,
$$
where the kernel $K$ satisfies the size and smoothness estimates
$$
|K(x,y)|
\le
\frac{C}{|x-y|^n}, \qquad x\neq y,
$$
and
$$
|K(x,y)-K(x',y)|+ |K(y,x)-K(y,x')|
\le
C\frac{|x-x'|^\delta}{|x-y|^{n+\delta}},
$$
for  all $|x-y|>2|x-x'|$.

Associated with $T$ is the truncated maximal operator
$$
T_\star f(x)
=
\sup_{0<\epsilon<\epsilon'<\infty}\Big| \int_{\epsilon<|x-y|<\epsilon'} K(x,y)f(y)dy \Big|.
$$

Let $M$ denote the Hardy-Littlewood maximal operator, that is,
$$
Mf(x)
=
\sup_{Q\ni x}\avgint_Q |f(y)|dy
=
\sup_{Q\ni x}\frac1{|Q|}\int_Q |f(y)|dy.
$$
where the supremum is taken over all cubes in $\re^n$ with sides parallel to the coordinate axes.

\medskip

\begin{theorem} \label{thm:main-thm:CZO}
Given a Calder\'on-Zygmund operator $T$, let $1<p<q<\infty$ and let $(u,v)$ be a pair of weights. If
the maximal operator satisfies
\begin{equation} \label{M:2-weight}
M: L^p(v) \rightarrow L^q(u) \quad \text{ and } \quad
M: L^{q'}(u^{1-q'}) \rightarrow L^{p'}(v^{1-p'}),
\end{equation}
then
\begin{equation} \label{T:2-weight}
\|Tf \|_{L^q(u)} \leq C\|f\|_{L^p(v)}
\quad \text{ and } \quad
\|T_\star f\|_{L^q(u)} \leq C\|f\|_{L^p(v)}.
\end{equation}
Analogously,
if
the maximal operator satisfies
\begin{equation} \label{M:2-weight:weak}
M: L^{q'}(u^{1-q'}) \rightarrow L^{p'}(v^{1-p'}),
\end{equation}
then
\begin{equation} \label{T:2-weight:weak}
\|Tf \|_{L^{q,\infty}(u)} \leq C\|f\|_{L^p(v)}
\quad \text{ and } \quad
\|T_\star f\|_{L^{q,\infty}(u)} \leq C\|f\|_{L^p(v)}.
\end{equation}

\end{theorem}

\bigskip

If the pairs of weights $(u,v)$ satisfy  any of the conditions in \eqref{M:2-weight}, then the weights $u$ and $v^{1-p'}$ are locally integrable.  This is a consequence of a characterization of the two weight norm inequalities for the maximal operator due to Sawyer~\cite{sawyer82b}.  He proved that the $L^p -L^q$ inequality holds if and only if for every cube $Q$,
\[ \left(\int_Q M(v^{1-p'}\chi_Q)(x)^q u(x)\,dx\right)^{1/q} \leq C\left(\int_Q v(x)^{1-p'}\,dx\right)^{1/p}<\infty, \]
and the $L^{q'}-L^{p'}$ inequality holds if and only if
\[ \left(\int_Q M(u\chi_Q)(x)^{p'}
  v(x)^{1-p'}\,dx\right)^{1/p'}
 \leq C\left(\int_Q u(x)\,dx\right)^{1/q'}<\infty.  \]
It is straightforward to construct pairs of weights that satisfy these
conditions.  For instance, in $\re$ both of these conditions follow easily for every $1<p\le q<\infty$
and the pair of weights $(u,v)$ with $u=\chi_{[0,1]}$ and $v^{-1}=\chi_{[2,3]}$ (i.e., $v=1$ in $[2,3]$ and $v=\infty$ elsewhere). Indeed, we only need to check Sawyer's inequalities for cubes Q that intersect both $[0,1]$ and $[2,3]$, in which case we have
$M(\chi_{[2,3]\cap Q})(x)\le |[2,3]\cap Q|$ for every $x\in [0,1]\cap Q$, and $M(\chi_{[0,1]\cap Q})(x)\le |[0,1]\cap Q|$ for every $x\in [2,3]\cap Q$. These readily imply the desired estimates.

\subsection*{Dyadic Calder\'on-Zygmund operators}
 A generalized dyadic grid  $\D$ in $\re^n$ is a set of generalized
 dyadic cubes with the following properties: if $Q\in \D$
then $\ell(Q)=2^{k}$, $k\in \Z$; if $Q, R\in \D$ and $Q\cap
R\neq\emptyset$ then $Q\subset R $ or $R\subset Q$; the cubes in $\D$
with $\ell(Q)=2^{-k}$ form a disjoint partition of $\re^n$ (see
\cite{lernerP2011-Dya} and \cite{lernerP2011-A2} for more details).

We say that $g_Q$ is a generalized a Haar function associated with $Q\in\D$ if
\begin{enumerate}\itemsep0.2cm
\item $\supp(g_Q) \subset Q$;
\item if $Q'\in \D$ and $Q'\subsetneq Q$, then $g_Q$ is constant on $Q'$;
\item $\|g_Q\|_\infty \leq 1$.
\end{enumerate}

Given a dyadic grid $\D$ and a pair $(m,k)\in \Z^2_+$, a linear operator $\S$ is
a generalized Haar shift operator (that is,  a dyadic Calder\'on-Zygmund
operator) of complexity type $(m,k)$ if it is bounded on $L^2(\re^n)$ and
$$
\S f(x)= \sum_{Q\in \D} \S_Q f(x)=
\sum_{Q\in \D}
\sum_{\substack{Q'\in\D_m(Q)\\ Q''\in \D_k(Q)}}
\frac{\langle f, g_{Q'}^{Q''} \rangle}{|Q|} g_{Q''}^{Q'}(x),
$$
where $\D_j(Q)$ stands for the dyadic subcubes of $Q$ with side length
$2^{-j}\ell(Q)$, $g_{Q'}^{Q''}$ is a generalized a Haar function
associated with $Q'$ and $g_{Q''}^{Q'}$ is a generalized a Haar function
associated with $Q''$. We say that the complexity of $\S$ is
$\kappa=\max(m,k)$. We also define the truncated Haar shift operator
$$
\S_\star f(x)=\sup_{0<\epsilon<\epsilon'<\infty} |\S_{\epsilon,\epsilon'} f(x)|
=
\sup_{0<\epsilon<\epsilon'<\infty} \Big|\sum_{\substack{Q\in\D\\\epsilon\le \ell(Q)\le \epsilon'}} \S_{Q} f(x)\Big|.
$$

An important example of a Haar shift operator on the real line is the Haar
shift (also known as the dyadic Hilbert transform)  $H^d$, defined
by
\[ H^df(x) = \sum_{I\in \Delta} \langle f, h_I\rangle \big(
h_{I_-}(x) - h_{I_+}(x)\big), \]
where, given a dyadic interval $I$, $I_+$ and $I_-$ are
its right and left halves, and
\[ h_I(x)  =  |I|^{-1/2}\big(\chi_{I_-}(x) - \chi_{I_+}(x)\big). \]
After renormalizing,  $h_I$ is a Haar function on $I$ and one can write $H^d$ as a
generalized Haar shift operator of complexity $1$.  These operators
have played a very important role in the proof of the $A_2$
conjecture:
see~\cite{DCU-JMM-CP2012,hytonenP2010,hytonen-perez-treil-volbergP}
and the references they contain for more information.

Associated with the dyadic grid $\D$ is the dyadic maximal function
$$
M_\D f(x)=
\sup_{x\in Q\in\D} \avgint_Q |f(y)|dy.
$$
Note that $M_\D$ is dominated pointwise by the Hardy-Littlewood
maximal operator.

\medskip

We can now state our result for dyadic Calder\'on-Zygmund
operators.

\begin{theorem} \label{thm:main-thm:dya}
Let $\S$ be a generalized Haar shift operator of complexity $\kappa$. Given $1<p<q<\infty$ and a pair of weights $(u,v)$, if
the dyadic maximal operator satisfies
\begin{equation} \label{M:2-weight:dya}
M_\D : L^p(v) \rightarrow L^q(u) \quad \text{ and } \quad
M_\D: L^{q'}(u^{1-q'}) \rightarrow L^{p'}(v^{1-p'}),
\end{equation}
then
\begin{equation} \label{T:2-weight:dya}
\|\S f\|_{L^q(u)} \leq C\kappa^2 \|f\|_{L^p(v)}
\quad \text{ and } \quad
\|\S_\star f\|_{L^q(u)} \leq C\kappa^2\|f\|_{L^p(v)}.
\end{equation}
Analogously,
if
the dyadic maximal operator satisfies
\begin{equation} \label{M:2-weight:weak:dya}
M_\D: L^{q'}(u^{1-q'}) \rightarrow L^{p'}(v^{1-p'})
\end{equation}
then
\begin{equation} \label{T:2-weight:weak:dya}
\|\S  f \|_{L^{q,\infty}(u)} \leq C \kappa^2 \|f\|_{L^p(v)}
\quad \text{ and } \quad
\|S_\star f\|_{L^{q,\infty}(u)} \leq C\kappa^2\|f\|_{L^p(v)}.
\end{equation}
\end{theorem}

\section{Proofs of the Main results}

\subsection*{Proof of Theorem \ref{thm:main-thm:CZO}}

We will prove our estimates for $T_\star$; the ones for $T$ are completely analogous.

Given a dyadic grid $\D$ we say that $\{Q_{j}^k\}_{j,k}$ is a
\textit{sparse family} of dyadic cubes if for any $k$ the cubes $\{Q_j^k\}_j$ are
pairwise disjoint; if $\Omega_k:=\cup_j Q_j^k$, then
$\Omega_{k+1}\subset \Omega_k$; and $|\Omega_{k+1}\cap Q_{j,k}|\le
\frac12 |Q_j^k|$. Given $\D$ and a sparse family
$\mathscr{S}=\{Q_j^k\}_{j,k}\subset\D$, define the positive dyadic
operator $\mathscr{A}$ by
$$
\mathscr{A} f(x)
=
\mathscr{A}_{\D, \mathscr{S}} f(x)
=
\sum_{j,k} f_{Q_j^k}\chi_{Q_j^k}(x)
$$
where $f_Q=\avgint_Q f(y)dy$.

For our proof we will use the main result in \cite{lernerP2011-Dya,lernerP2011-A2}.   Given a Banach function space $X$ and a
non-negative function $f$,
\begin{equation}\label{main-Lerner}
\|T_\star f\|_X
\le
C(T,n)\sup_{\D,\mathscr{S}} \|\mathscr{A}_{\D, \mathscr{S}} f\|_X,
\end{equation}
where the supremum is taken over all dyadic grids $\D$ and sparse families
$\mathscr{S}\subset \D$.  To prove Theorem~\ref{thm:main-thm:CZO} we
apply this result with $X=L^q(u)$ or $X=L^{q,\infty}(u)$; it will then
suffice to show that our assumptions on $M$ guarantee that
$\mathscr{A}_{\D, \mathscr{S}}$ satisfies the corresponding two weight
inequalities.

To prove this fact we will use a result by Lacey, Sawyer and
Uriate-Tuero~\cite{LSUT-P2010}.  Given
a sequence of non-negative constants $\alpha=\{\alpha_Q\}_{Q\in\D}$, define the positive operator
$$
T_\alpha f(x)
=
\sum_{Q\in\D} \alpha_Q f_Q\chi_Q(x).
$$
Further, given $R\in\D$ we define the ``outer truncated'' operator
$$
T_{\alpha}^{R}f(x)
=
\sum_{\substack{Q\in\D \\Q\supset R}} \alpha_Q f_Q\chi_Q(x).
$$
In \cite{LSUT-P2010} it was shown that for all $1<p<q<\infty$,
$T_\alpha:L^p(v)\rightarrow L^{q}(u)$ if and only if there exist constants
$C_1$ and $C_2$ such that for every $R\in \D$
\begin{equation}\label{test:TR}
\left(\int_{\re^n} T_\alpha^R (v^{1-p'}\chi_R)(x)^qu(x)dx\right)^{\frac1q}
\le
C_1
\left(\int_{R} v(x)^{1-p'}dx\right)^{\frac1p},
\end{equation}
and
\begin{equation}\label{test:TR:dual}
\left(\int_{\re^n} T_\alpha^R (u\chi_R)(x)^{p'}v(x)^{1-p'}dx\right)^{\frac1{p'}}
\le
C_2
\left(\int_{R} u(x)dx\right)^{\frac1{q'}}.
\end{equation}
Furthermore, for $1<p<q<\infty$, $T_\alpha:L^p(v)\rightarrow L^{q,\infty}(u)$  holds
if and only if there exists a constant $C_2$ such that for every $R\in \D$, \eqref{test:TR:dual} holds.

We can apply these results to the operator
$\mathscr{A}=\mathscr{A}_{\D,\mathscr{S}}$ where $\D$ and
$\mathscr{S}$ are fixed, since $\mathscr{A}=T_\alpha$ with
$\alpha_Q=1$ if $Q\in\mathscr{S}$ and $\alpha_Q=0$ otherwise. Fix $R\in\D$; to
estimate $\mathscr{A}^R$, take the increasing family of cubes
$R=R_0\subsetneq R_1\subsetneq R_2\subsetneq \dots$ with $R_k\in\D$
and $\ell(R_k)=2^k \ell(R)$.  Define $R_{-1}=\emptyset$.  Note that
$\supp \mathscr{A}^R\subset \cup_{k\ge 0} R_k$. Then for every
non-negative function $f$ and for every $x\in R_k\setminus R_{k-1}$
with $k\ge 0$ we have that
\begin{multline*}
0
\le \mathscr{A}^R (f\chi_R)(x)
\le
\sum_{j=0}^\infty (f\chi_R)_{R_j}\chi_{R_{j}}(x)
=
f_R
\sum_{j=k}^\infty 2^{-j\,n}
\\
\lesssim
f_R 2^{-k\,n}
=
(f\chi_R)_{R_k}
\le
M_\D (f\chi_R)(x).
\end{multline*}
Consequently, for every $x\in \re^n$,
\begin{equation}\label{A-M}
0
\le \mathscr{A}^R (f\chi_R)(x)
\lesssim
M_\D (f\chi_R)(x) \leq M(f\chi_R)(x).
\end{equation}

Inequality \eqref{A-M} together with our hypothesis~\eqref{M:2-weight}
implies \eqref{test:TR} and \eqref{test:TR:dual}.  Therefore, we have
that $\mathscr{A}:L^p(v)\rightarrow L^{q}(u)$ with constants depending
on the dimension, $p$, $q$ and the implicit constants in
\eqref{M:2-weight}. Therefore, by Lerner's estimate
\eqref{main-Lerner} we get $T_\star:L^p(v)\rightarrow L^{q}(u)$
as desired.

For the weak-type estimates we proceed in the same manner, using the
fact that \eqref{M:2-weight:weak} yields \eqref{test:TR:dual} and
therefore $\mathscr{A}:L^p(v)\rightarrow L^{q,\infty}(u)$. This in
turn implies, by Lerner's estimate \eqref{main-Lerner} applied to
$X=L^{q,\infty}(u)$, that $T_\star:L^p(v)\rightarrow
L^{q,\infty}(u)$.

\subsection*{Proof of Theorem \ref{thm:main-thm:dya}}

Fix $\D$ and a generalized Haar shift operator of complexity
$\kappa$. As before we can work with $\S_\star$.   We can repeat the
previous argument except that we want to keep the fixed dyadic structure $\D$. A careful
examination of \cite[Section 5]{lernerP2011-Dya} shows that, given a
Banach function space $X$, we have
\begin{equation}\label{main-Lerner:dya}
\|\S_\star f\|_X
\le
C_n\kappa^2\sup_{\mathscr{S}} \|\mathscr{A}_{\D, \mathscr{S}} f\|_X,
\qquad f\ge 0,
\end{equation}
where the supremum is taken over all sparse families $\mathscr{S}\subset \D$. We
emphasize that in \cite[Section 5]{lernerP2011-Dya} there is an
additional supremum over the dyadic grids $\D$. This is because at some
places the dyadic maximal operator is majorized by the regular
Hardy-Littlewood maximal operator and the latter is in turn controlled
by a sum of $\mathscr{A}_{\D_\alpha, \mathscr{S}_\alpha}$ for $2^n$
dyadic grids $\D_\alpha$. However, keeping $M_\D$ one can easily show
that \eqref{main-Lerner:dya} holds. Details are left to the interested
reader.

Given \eqref{main-Lerner:dya}, we fix a sparse family
$\mathscr{S}\subset \D$ and write $\mathscr{A}=\mathscr{A}_{\D,
  \mathscr{S}}$. Arguing exactly as before we obtain \eqref{A-M}.
Thus, \eqref{M:2-weight:dya} implies \eqref{test:TR} and
\eqref{test:TR:dual} and therefore the result from \cite{LSUT-P2010} yields
$\mathscr{A}:L^p(v)\rightarrow L^{q}(u)$ with constants depending
on the dimension, $p$, $q$ and the implicit constants in
\eqref{M:2-weight:dya}. Combining this with Lerner's estimate
\eqref{main-Lerner:dya} applied to $X=L^q(u)$ we conclude as desired
that $\S_\star :L^p(v)\rightarrow L^{q}(u)$.  We get the
weak-type estimate by adapting the above proof in exactly the same
way.

\bibliographystyle{plain}
\bibliography{testing}

\begin{thebibliography}{10}

\bibitem{MR2351373}
D.~Cruz-Uribe, J.~M. Martell, and C.~P{\'{e}}rez.
\newblock Sharp two-weight inequalities for singular integrals, with
  applications to the {H}ilbert transform and the {S}arason conjecture.
\newblock {\em Adv. Math.}, 216(2):647--676, 2007.

\bibitem{MR2628851}
D.~Cruz-Uribe, J.~M. Martell, and C.~P{\'{e}}rez.
\newblock Sharp weighted estimates for approximating dyadic operators.
\newblock {\em Electron. Res. Announc. Math. Sci.}, 17:12--19, 2010.

\bibitem{cruz-martell-perezBook}
D.~Cruz-Uribe, J.M. Martell, and C.~P{\'{e}}rez.
\newblock {\em Weights, extrapolation and the theory of {R}ubio de {F}rancia},
  volume 215 of {\em Operator Theory: Advances and Applications}.
\newblock Birkh\"auser/Springer Basel AG, Basel, 2011.

\bibitem{DCU-JMM-CP2012}
D.~Cruz-Uribe, J.M. Martell, and C.~P\'erez.
\newblock Sharp weighted estimates for classical operators.
\newblock {\em Adv. in Math}, 229:408--441, 2012.

\bibitem{cruz-reznikov-volberg}
D.~Cruz-Uribe, A.~Reznikov, and A.~Volberg.
\newblock Logarithmic bump conditions and the two-weight boundedness of
  {C}alder\'on--{Z}ygmund operators.
\newblock {\em Preprint}, 2012.
\newblock arXiv:1112.0676.

\bibitem{hytonenP2010}
T.~Hyt\"{o}nen.
\newblock The sharp weighted bound for general {C}alder\'{o}n-{Z}ygmund
  operators.
\newblock {\em Ann. of Math. (2)}, (to appear).
\newblock arXiv:1007.4330 (2010).

\bibitem{hytonen-perez-treil-volbergP}
T.~Hyt\"{o}nen, C.~P\'{e}rez, S.~Treil, and A.~Volberg.
\newblock Sharp weighted estimates of the dyadic shifts and ${A}_2$ conjecture.
\newblock {\em Preprint}, 2010.
\newblock arXiv:1010.0755.

\bibitem{LSUT-P2010}
M.~Lacey, E.~Sawyer, and I.~Uriarte-Tuero.
\newblock Two weight inequalities for discrete positive operators.
\newblock {\em Preprint}, 2010.
\newblock arXiv:0911.3437.

\bibitem{lernerP2011-Dya}
A.~Lerner.
\newblock On an estimate of {C}alder\'{o}n-{Z}ygmund operators by dyadic
  positive operators.
\newblock {\em Preprint}, 2012.
\newblock arXiv:1202.1860.

\bibitem{lernerP2011-A2}
A.~Lerner.
\newblock A simple proof of the ${A}_2$ conjecture.
\newblock {\em Preprint}, 2012.
\newblock arXiv:1202.2824.

\bibitem{reguera-scurry}
M.~Reguera and J.~Scurry.
\newblock On joint estimates for maximal functions and singular integrals in
  weighted spaces.
\newblock {\em Preprint}, 2011.
\newblock arXiv:1109.2027.

\bibitem{sawyer82b}
E.T. Sawyer.
\newblock A characterization of a two-weight norm inequality for maximal
  operators.
\newblock {\em Studia Math.}, 75(1):1--11, 1982.

\end{thebibliography}

\end{document}